\documentclass[a4paper,11pt]{article}
\usepackage[utf8]{inputenc}

\usepackage[left=1.8cm, right=1.8cm, top=2.2cm, bottom=2.5cm]{geometry}
\usepackage{amsthm}

\usepackage{comment}

\usepackage[
backend=biber,
style=alphabetic,
sorting=anyt
]{biblatex}
\usepackage{dirtytalk}
\usepackage{amsmath}
\usepackage{amssymb}
\usepackage{enumitem}
\usepackage{verbatim}
\usepackage{stmaryrd}
\usepackage{amsfonts}
\usepackage{mathtools}
\usepackage{titlesec}
\usepackage{authblk}
\usepackage[T2A]{fontenc}
\usepackage[english]{babel}

\DeclareFieldFormat{pages}{#1}

\titlespacing\section{0pt}{12pt plus 4pt minus 2pt}{2pt plus 2pt minus 2pt}

\makeatletter
\renewrobustcmd*{\mkbibemph}{}
\protected\long\def\blx@imc@mkbibemph#1{#1}
\makeatother

\newtheoremstyle{boldremark}
    {\dimexpr\topsep/2\relax} 
    {\dimexpr\topsep/2\relax} 
    {}          
    {}          
    {\bfseries} 
    {.}         
    {.5em}      
    {}          

\theoremstyle{boldremark}

\newtheorem{Lem}{Lemma}
\newtheorem*{AppLem}{Approximation lemma}
\newtheorem*{FinTh}{Finite-normal-subgroup theorem}
\newtheorem*{RetrTh}{Nilpotent-strong-retract theorem}

\newtheorem{Cons}{Corollary}
\newtheorem*{Cons2}{Corollary 2 \cite{verbal}}

\newtheorem*{Proof}{Proof}
\newtheorem{Question}{Question}
\newtheorem{Prop}{Proposition}

\DeclareMathOperator{\ord}{ord}

\DeclareMathOperator{\Aut}{Aut}

\DeclareMathOperator{\var}{\mathbf{var}}
\DeclareMathOperator{\End}{End}

\let\leq\leqslant
\let\geq\geqslant
\let\phi\varphi
\let\epsilon\varepsilon
\let\tilde\widetilde

\title{\textsc{\textbf{\large{finite normal subgroups of strongly verbally closed groups\vspace{-0.5cm}}}}}

\author{Filipp D. Denissov \\ \textsl{Faculty of Mathematics and Mechanics of Moscow State University} \\ \textsl{Moscow 119991, Leninskie gory, MSU.} \\ \textsl{Moscow center for Fundamental and Applied Mathematics.} \\
\textsl{denissov.filipp@gmail.com}}
\date{}

\begin{document}

\renewcommand{\abstractname}{\vspace{-\baselineskip}}

\maketitle

\vspace{-1.5cm}

\begin{abstract}
In the recent paper by A. A. Klyachko, V. Yu. Miroshnichenko, and A. Yu. Olshanskii, it is proven that the center of any finite strongly verbally closed group is its direct factor. One of the results of the current paper is the generalization of this nontrivial fact to the case of finite normal subgroups of any strongly verbally closed groups. It follows from this generalization that finitely generated nilpotent groups with nonabelian torsion subgroups are not strongly verbally closed.
\end{abstract}

\section*{\normalsize 1. Introduction}

\noindent A subgroup $H$ of a group $G$ is called \textit{verbally closed} \cite{myasnikov} if any equation of the form
\begin{equation*}
    w(x_1, x_2, \dots, x_n) = h, \text{ where } w  \text{ is an element of the free group } F(x_1, \dots, x_n) \text{ and } h \in H,
\end{equation*}
having solutions in $G$ has a solution in $H$. If each system of equations with coefficients from H
\begin{equation*}
    \{w_1 (x_1, \dots) = 1, \dots, w_m (x_1, \dots) = 1 \}, \text{ where } w_i \in H \ * \ F(x_1, \dots, x_n) \text{ (and } * \text{ means the free product),}
\end{equation*}
having solutions in $G$ has a solution in $H$, then the subgroup $H$ is called \textit{algebraically closed} in $G$. Note that if the subgroup $H$ is algebraically closed in the group $G$, then it is verbally closed in $G$.

A group $G$ is called \textit{strongly verbally closed} if it is algebraically closed in any group containing $G$ as a verbally closed subgroup. Thus, the verbal closedness (as well as the algebraic closedness) is a property of a subgroup, while the strong verbal closedness is a property of an abstract group. The class of strongly verbally closed groups is fairly wide. For example, it includes \begin{itemize}[topsep=0pt]
\itemsep-0.5em
\renewcommand\labelitemi{---}
\item all abelian groups \cite{mazhuga},
\item all free groups \cite{free},
\item all virtually free groups containing no nontrivial finite normal subgroups \cite{free}, \cite{virtually},
\item all groups decomposing nontrivially into a free product \cite{mazhuga19},
\item fundamental groups of all connected surfaces except the Klein bottle \cite{mazhuga}, \cite{klein},
\item all finite groups with nonabelian monolith \cite{verbal},
\item inifinite dihedral group \cite{virtually} and any finite dihedral group whose order is not divisible by $8$ \cite{verbal},
\item all acylindrically hyperbolic groups with no nontrivial finite normal subgroups \cite{bog18}.
\end{itemize}
The class of non-strongly-verbally-closed groups is fairly wide too. Among such groups are the following:
\begin{itemize}[topsep=0pt]
\itemsep-0.5em
\renewcommand\labelitemi{---}
\item the already mentioned fundamental group of the Klein bottle \cite{klein},
\item the discrete Heisenberg group \cite{verbal},
\item any finite group, whose center is not its direct factor (in particular, any finite nonabelian nilpotent group) \cite{verbal}, \cite{nilpotent}, \cite{free}.
\end{itemize}
Proving the strong verbal closedness (as well as its absence) of a group is not easy. In \cite{verbal}, for example, a question is raised:

\begin{Question}
\textsl{Does there exist a finitely generated nilpotent nonabelian strongly verbally closed group?}
\end{Question}

A negative answer to this question would yield a broad generalization of the last two examples of non-strongly-verbally-closed groups mentioned above. So far, we managed to give a partial answer to this question. More precisely, we proved the absence of strong verbal closedness of finitely generated nilpotent groups with nonabelian torsion subgroups and of some finitely generated nilpotent nonabelian groups with abelian torsion subgroups.

A property that is stronger than the strong verbal closedness is the property of being a strong retract \cite{verbal}. A group $H$ is called a \textit{strong retract} if it is a retract of any group $G \geq H$ from the variety generated by the group $H$.

\noindent Let us recall some terminology \cite{variety}:
\begin{itemize}[topsep=0pt]
\itemsep-0.5em
\renewcommand\labelitemi{---}
\item \textsl{the }\textit{variety generated by a class of groups} $\mathcal{K}$ is the class of all groups satisfying all identities that hold in all groups from $\mathcal{K}$,
\item the variety generated by a group $G$ is designated by $\var G$.
\end{itemize}
This gives rise to the following question from \cite{verbal}:

\begin{Question}
\textsl{What is an arbitrary finite strong retract?}
\end{Question}

In \cite{verbal} some examples of strong retracts are provided. In the next section, we describe all the nilpotent strong retracts.

\vspace{1mm}

\noindent Below we provide a brief list of \textbf{notation} we use.

If $x,y$ are elements of some group, then the symbol $[x,y]$ denotes their commutator $x^{-1} y^{-1} x y$. The symbol $\ord (x)$ denotes the order of an element $x$ of a group $G$. The center of a group $G$ is denoted by $Z(G)$, and its commutator subgroup is denoted by $G'$. The centralizer of a subset $X$ of a group $G$ is denoted by $C(X)$. The symbol $\langle \langle X \rangle \rangle$ stands for the normal closure of a subset $X$ of a group $G$ (that is the intersection of all normal subgroups of $G$ containing $X$). The free group with a basis $X$ is denoted as $F(X)$ or $F_n$ in case $X$ has $n \in \mathbb{N}$ elements. Identical mapping from $X$ to itself is denoted by $id$. We use the symbol $H \cong G$ to express the fact that groups $H$ and $G$ are isomorphic. Finally, the symbol $H \leq G$ denotes the fact that a group $H$ is a subgroup of $G$. The symbol $H \unlhd G$ denotes the fact that $H$ is a normal subgroup of $G$.

\vspace{1mm}

The author is grateful to his supervisor Anton Alexandrovich Klyachko for formulation of the problem and for valuable remarks during the work.

\section*{\normalsize 2. Nilpotent strong retracts}

\noindent Note that in case when $G$ is an abelian group, $H \leq G$ is its retract if and only if $H$ is a direct summand of $G$. It means that the property of being a strong retract for the abelian group $G$ is equivalent to the property of $G$ being a direct summand of any group $H \in \var G$ containing $G$. For the further discussion, we need the description of all varieties of abelian groups (see \cite{fuks}, paragpaph $18$, exercise $7$):
\begin{quote}
    \textsl{Varieties of abelian groups are precisely the following classes of groups: 1) the class of all abelian groups; 2) the class of all abelian groups of a period divising $n \in \mathbb{N}$.}
\end{quote} 

\noindent \noindent 
Recall that the \textit{period of a group $G$} is the least number $n \in \mathbb{N}$, such that $x^n = 1$ for any $x \in G$. If such a number exists, then $G$ is a group of \textit{bounded period}.

To begin with, consider the case, when $G$ is not a group of bounded period. Then, according to the description, $\var G$ is the class of all abelian groups. The following is true of divisible abelian groups (see, for example, \cite{kurosh}):
\begin{quote}
    \textsl{If $G$ is a divisible abelian group, and $H$ is an abelian group such that $G \leq H$, then $G$ is a direct summand of $H$.}
\end{quote}

\noindent Let us remind that a group $G$ is called \textit{divisible} if for any $ g \in G$ and $n \in \mathbb{N}$, the equation $x^n = g$ has a solution in $G$. 

\begin{Prop}
\textsl{An abelian group $G$ of unbounded period is a strong retract if and only if it is divisible.}
\end{Prop}

\begin{Proof}
Sufficiency follows from the fact provided above. Let $G$ be an abelian group of unbounded period. Then, as it was noted earlier, $\var G$ is the class of all abelian groups. In particular, $\var G$ contains a divisible group $H$ containing $G$ \cite{kurosh}. Though, if $G$ is not divisible itself, it is not a direct summand of $H$ (as direct summands of a divisible group are divisible themselves \cite{kurosh}), so $G$ is not a strong retract. \qedsymbol
\end{Proof}

\noindent Let us move on to abelian groups of bounded period. The first Prüfer theorem provides a complete description of these groups \cite{kurosh}:
\begin{quote}
    \textsl{An abelian group $G$ of bounded period $d$ is a direct sum of primary cyclic groups, i.e. $G \cong \bigoplus_{i \in I} \mathbb{Z}_{p_{i}^{k_i}}$, where $p_i$ are prime numbers and $k_i$ are natural numbers such that $p_i^{k_i} \vert d$, $i \in I$ ($I$ is an index set).}
\end{quote}

\noindent We need the following variation of the Zorn's lemma \cite{fuks}:

\begin{quote}
    \textsl{Let $M \neq \varnothing$ be a partially ordered set. Suppose that every chain in $M$ (a totally ordered subset of $M$) has an upper bound. Then $M$ contains a maximal element.}
\end{quote}

\noindent Now, we are ready to proceed with our description:

\begin{Prop}
\textsl{An abelian group $G$ of bounded period is a strong retract if and only if in its decomposition into the direct sum of primary cyclic groups, orders of any distinct direct summands are either equal or coprime:}
\begin{equation*}
    G \cong \bigoplus_{i = 1}^m C_{p_i^{k_i}} (n_i), \text{ where } C_{p_i^{k_i}} (n_i) \text{ is equal to the direct sum of } n_i \text{ copies of the group } \mathbb{Z}_{p_i^{k_i}},
\end{equation*}

\noindent \textsl{where all prime numbers $p_i$ are distinct, $m, k_i \in \mathbb{N}$, and $n_i$ are some cardinal numbers.}
\end{Prop}

\begin{Proof}
Suppose that $G$ cannot be decomposed into such a direct sum. We may assume that

\begin{equation}
    G = \bigoplus_{i = 1}^m \bigoplus_{j \in I_i} \mathbb{Z}_{p_i^{k_j}},
\end{equation}

\noindent where $m \in \mathbb{N}$, $|I_i| = n_i$ and among $k_j$, $j \in I_i$ there are only finitely many different ones (because $G$ is a group of bounded period) but there exists $i \in \{1, \dots, m\}$ such that for some $j_1, j_2 \in I_i,$ $k_{j_1} \neq k_{j_2}$. 

Consider the group: $\displaystyle H = \bigoplus_{i = 1}^m C_{p_i^{s_i}} (n_i)$, where 
\begin{equation*}
    s_i = \max \{k_j \ \vert \ \mathbb{Z}_{p_i^{k_j}} \text{ is a direct summand in the decomposition (1)}\}, \ i = 1,2,\dots, m.
\end{equation*}
Since both $G$ and $H$ are of the same period $\prod_{i = 1}^m s_i$, it follows from the description of abelian varieties that $H \in \var G$ .

Consider the injection $f : G \to H$, which works on each direct summand from $(1)$ as follows: let $i \in \{1, \dots, m\}$, $j \in I_i$, $f : \mathbb{Z}_{p_i^{k_j}} \hookrightarrow \mathbb{Z}_{p_i^{s_i}}$, where $\mathbb{Z}_{p_i^{s_i}}$ is the $j$th summand from the decomposition of $C_{p_i^{s_i}} (n_i)$ into the direct sum. Every direct summand from $(1)$ is mapped into the corresponding direct summand of the decomposition of $H$, so that the restriction of $f$ to $\mathbb{Z}_{p_i^{k_j}}$ is a natural injection: if $k_j = s_i$, then it is the identical map; otherwise it is a mapping to the subgroup of $\mathbb{Z}_{p_i^{s_i}}$ of the order $p_i^{k_j}$. From the uniqueness of the decomposition of an abelian group of bounded period into the direct sum of primary cyclic groups \cite{fuks}, it follows that $f(G)$ is not a direct summand of $H$. Thus, $G$ is not a strong retract.

Now, suppose that $G$ has the decomposition from the statement of the theorem. Let $H \in \var G$ and let $f : G \hookrightarrow H$  be a monomorphism. As any monomorphism preserves the order of an element, the $p_i$th component of $G$ is mapped into the $p_i$th component of $H$ under $f$, so it suffices to prove the theorem only for the case $G = C_{p^{k}} (n)$, where $p$ is prime, $k \in \mathbb{N}$, and $n$ is some cardinal number. 

Let us show that there exists such $X \leq H$ that $H = f(G) \oplus X$. In Zorn's lemma, take the set of all subgroups of $H$ having trivial intersection with $f(G)$ as $M$:
\begin{equation*}
    M = \{Y \leq H \ \vert \ Y \cap f(G) = \{0\}\}.
\end{equation*}

\noindent Order on $M$ is introduced as follows: for $X, Y \in M$, $X \leqslant Y$ if $X$ is a subgroup of $Y$. It can be verified directly that this is an order on $M$. Set $M$ is nonempty: $\{0\} \in M$. Any chain $\{Y_{\alpha}\} \subseteq M$ of subgroups having trivial intersection with $f(G)$ is bounded by an element $Y \in M$, where $Y = \cup_{\alpha} Y_{\alpha}$. Consequently, Zorn's lemma is applicable, and $M$ contains a maximal element $X$: $X \leq H$, $X \cap f(G) = \{0\}$, and $X$ is not a subgroup of any bigger (relatively to the order we introduced) subgroup satisfying this property.

From $X \cap f(G) = \{0\}$ it follows that $f(G) + X = f(G) \oplus X$. It remains to prove that $H = f(G) + X$. Let $h \in H$. There exists such $k \in \mathbb{N}$ that $kh \in f(G) + X$. Indeed, otherwise $\langle h \rangle \cap (f(G) + X) = \{0\}$, which means that $(\langle h \rangle + X) \cap f(G) = \{0\}$, leading to a contradiction with the maximality of $X$. 

Let $s$ be the least of such numbers $k$. Without loss of generality, assume that $s$ is prime or that $s = 1$ (otherwise, take a power of $h$ instead of $h$). Two cases are possible:

1) $s = p$. Then, $ph = f(g) + x$ for some $g \in G$, $x \in X$. If $g = p g_1$, $g_1 \in G$  ($g_1$ may be equal to zero), then $ph - f(p g_1) = x$. However, from $h - f(g_1) \not\in X$ (as $h \not \in f(G) + X$) it can be obtained that $(X + \langle h - f(g_1) \rangle) \cap f(G) = \{0\}$, which leads to a contradiction with the maximality of $X$. Consequently, $g \neq p g_1$ for any $g_1 \in G$. As $g \neq 0$,  $\ord (g) = p^k$. Though, $\ord (ph) = p^r < p^k$, so $p^r (ph) = 0 = p^r(f(g)) + p^r x$. As the sum $f(G) + X$ is direct, $p^r f(g) = p^r x = 0$, which means that $p^r g = 0$, which is impossible.

2) $s \neq p$. For abelian groups of period $p$, the mapping $g \mapsto sg$ is an automorphism, so, as $sh = f(g) + x$ for some $g \in G$, $x \in X$, there exist such $g_1 \in G$, $x_1 \in X$ that $g = sg_1$, $x = s x_1$. Thus, $s(h - f(g_1) - x_1) = 0$. No nontrivial element of $H$ has the order of $s$, so $h = f(g_1) + x_1$. 

As a result, $H = f(G) \oplus X$, and $G$ is a strong retract. \qedsymbol

\end{Proof}

\begin{Prop}
\textsl{The center of a strong retract is its direct factor.}
\end{Prop}

\begin{Proof}
Let $G$ be a strong retract. The center of any group is a normal subgroup, so it suffices to prove that $Z(G)$ is a retract of $G$. Consider the central product of $G$ with its copy $\tilde{G}$ with joined center:
\begin{equation*}
    K = G \underset{Z(G) = Z(\Tilde{G})}{\times} \Tilde{G}= (G \times \Tilde{G}) / \{(g, g^{-1}) \vert g \in Z(G)\} \in \var G.
\end{equation*}
The group $\Tilde{G}$ is isomorphic to the group $G$, so it is a strong retract too. Let $\rho$ be a retraction of $K$ to its subgroup $\Tilde{G}$. From the fact that in the group $K$, the group $G$ commutes with the group $\Tilde{G}$, we obtain $\rho (G) \leq Z(G)$. By definition of the retraction, $\rho(g) = g$ is true for any element $g \in Z(G)$. Thus, the restriction of $\rho$ to the subgroup $G$ of the group $K$ is the desired retraction to $Z(G)$.  \qedsymbol
\end{Proof}

\noindent The following simple proposition shows that consideration of nilpotent groups does not yield any new strong retracts:

\begin{Prop}
\textsl{Nilpotent strong retract is an abelian group.}
\end{Prop}

\begin{Proof}
Any nontrivial normal subgroup of a nilpotent group intersects the center of this group nontrivially (see \cite{kargapolov}). From this fact and from the proposition $3$, we obtain that any nilpotent strong retract is equal to its center. \qedsymbol
\end{Proof}

\noindent As a result, we proved the following theorem:

\begin{RetrTh}
\textsl{Nilpotent strong retracts are precisely divisible abelian groups and abelian groups of bounded period in whose decomposition into the direct sum of primary cyclic groups, orders  of any distinct direct summands are either equal or
coprime.}
\end{RetrTh}

In the next paragraph we show that many nilpotent groups are not even strongly verbally closed.

\section*{\normalsize 3. Finite normal subgroups of strongly verbally closed groups}

\noindent We say that a group presentation $\langle X \ \vert \ R \rangle$ is \textit{finitely presented} over a group presentation $\langle Y \ \vert \ S \rangle$, if there exist such finite sets $A$ and $B$ that $\langle X \ \vert \ R \rangle \cong \langle X' \ \vert \ R' \rangle$, where $X' = Y \cup A$, $R' = S \cup B$.

The following lemma reveals that this definition is, in fact, a group property (which means it does not depend on the choice of a group presentation), so it makes sense to speak about the finite presentability of one group over the other group:

\begin{Lem}
\textsl{Suppose that a group presentation $\langle X \ \vert \ R \rangle$ is finitely presented over a group presentation $\langle Y \ \vert \ S \rangle$ and $\langle  Y \ \vert \ S \rangle \cong \langle Y' \ \vert \ S' \rangle$. Then $\langle X \ \vert \  R \rangle$ is finitely presented over $ \langle Y' \ \vert \ S' \rangle$.}
\end{Lem}

\begin{Proof}
We may assume that $X = Y \cup A$ and $R = S \cup B$ for some finite sets $A$ and $B$. It is known (see, for example, \cite{lindon}) that groups defined by group presentations $\langle Y  \ \vert \ S \rangle$ and $\langle Y' \ \vert \ S' \rangle$ are isomorphic if and only if presentation $\langle Y' \ \vert \ S' \rangle$ is obtained from presentation $\langle Y \ \vert \ S \rangle$ by applying a finite number of \textit{Tietze transformations}: \begin{itemize}[topsep=0pt]
\itemsep-0.5em
\renewcommand\labelitemi{---}
\item adding to the set $S$ an arbitrary set $T \subseteq \langle \langle S \rangle \rangle \unlhd F(Y)$ of its consequences,
\item adding to the set $Y$ an arbitrary set $\tilde{Y}$ while ading to $S$ a set $\{\tilde{y} = w_{\tilde{y}} \ \vert \ \tilde{y} \in \tilde{Y}, w_{\tilde{y}} \in F(Y)\}$,
\end{itemize}
and their inverses. It is sufficient to prove the lemma only for the case, when $\langle Y' \ \vert \ S' \rangle$ is obtained from $\langle Y \ \vert \ S \rangle$ by applying one Tietze transformation. One can easily verify that in case of the first transformation, $X' = X$ and $R' = R \cup T$, while in case of the second transformation, $X' = X \cup \tilde{Y}$ and $R' = R \cup \{\tilde{y} = w_{\tilde{y}} \ \vert \ \tilde{y} \in \tilde{Y}, w_{\tilde{y}} \in F(Y)\}$ provide the desired group presentation. \qedsymbol
\end{Proof}

\noindent By virtue of Lemma $1$, the following definition may be introduced:

A group $G$ is \textit{finitely presented} over a group $H$, if there exists such a presentation of $G$ that it is finitely presented over any presentation of $H$. 

\begin{Lem}
\textsl{Suppose that $G$ contains a subgroup $H$ and a finite normal subgroup $N$ such that $G / N$ is finitely presented over $H / (H \cap N)$. Then $G$ is finitely presented over $H$.}
\end{Lem}

\noindent \textbf{Proof} (with minor changes) replicates the proof of the Hall theorem \cite{hall} about preservation of finite presentability of a group under extensions (see also \cite{robinson}).

Let $G$ be a group, $H = \langle X \ \vert \ R \rangle \leq G$, and $N = \langle Y \ \vert \ S \rangle \ \unlhd \ G$ be its finite subgroup, where $Y$ and $S$ are finite sets. By condition of the lemma, the group $G / N$ is finitely presented over $H / (H \cap N) = \langle X \ \vert \ R \cup C \rangle$, where $\langle \langle C \rangle \rangle = H \cap N$ and the set $C$ is finite. Consequently
\begin{equation*}
    G / N \cong \langle X \cup A \ \vert \ R \cup C \cup B \rangle,
\end{equation*}
where sets $A$ and $B$ are finite. 

Let us construct a presentation of the group $G$. As the set of generators, take $\overline{X} \cup \overline{A} \cup \overline{Y}$, where sets $\overline{X}$, $\overline{A}$, $\overline{Y}$ are in one-to-one correspondence with sets $X$, $A$, $Y$ respectively. The sets $R$, $S$, $C$, and $B$ are in correspondence with the sets $\overline{R}$, $\overline{S}$, $\overline{C}$, and $\overline{B}$ respectively. As the set of defining relations, take the union of the following sets: $\overline{R}$, $\overline{S}$, $\overline{C}_1 = \{c w_c^{-1} \ \vert \ c \in \overline{C}, w_c \in F(\overline{Y})\}$,  $\overline{B}_1 = \{b w_b^{-1} \ \vert \ b \in \overline{B}, w_b \in F(\overline{Y})\}$ ($c \in \overline{C}$ and $b \in \overline{B}$ are considered as words from $F(\overline{X})$ and from $F(\overline{X} \cup \overline{A})$ respectively), $\overline{T} = \{a^{-1} y a w_{a,y}^{-1}, \ a y a^{-1} v_{a,y}^{-1} \ \vert \ a \in \overline{A}, y \in \overline{Y}, \  w_{a, y}, v_{a,y} \in F(\overline{Y})\}$:
\begin{equation*}
    \tilde{G} = \langle \overline{X} \cup \overline{A} \cup \overline{Y} \ \vert \  \overline{R} \cup \overline{S} \cup \overline{C}_1 \cup \overline{B}_1 \cup \overline{T} \rangle.
\end{equation*}
Consider a surjective homomorphism $\theta : \tilde{G} \to G$, defined with the following bijections $\overline{X} \to X$, $\overline{A} \to A$, $\overline{Y} \to Y$ on the generators (defining relations are mapped into true identities under such a map on generators, so such a homomorphism exists). The restriction $\theta \vert_{K} : K \to N$ on the subgroup $K = \langle \overline{Y} \rangle \leq \tilde{G}$ is an isomorphism as all the relations in the alphabet $\overline{Y}$ in $\tilde{G}$ are consequences of the defining relations $\overline{S}$. Besides, $K \unlhd \tilde{G}$. 

Homomorphism $\tilde{\theta} : \tilde{G} / K \to G / N$ generated by $\theta$, is an isomorphism too. Now, let $g \in \ker \theta$. Then $gK \in \ker \tilde{\theta}$, but $\tilde{\theta}$ is an isomorphism, so $g \in K$. Finally, $\theta \vert_{K}$ is an isomorphism, so $g = 1$. \qedsymbol

\vspace{1mm}

\noindent The following lemma provides a criterion for algebraic closedness of a subgroup $H$ of a group $G$ in case, when $G$ is finitely presented over $H$ (for similar propositions, refer to \cite{myasnikov}):

\begin{Lem}
\textsl{Suppose that $H = \langle X \ \vert \ R \rangle$ is a subgroup of $G$ and $G$ is finitely presented over $H$. The subgroup $H$ is algebraically closed in $G$ if and only if $H$ is a retract of $G$.}
\end{Lem}

\begin{Proof}
Suppose $H$ is algebraically closed in $G$ and $A = \{a_1, \dots, a_m\}$, $B = \{s_1, \dots, s_n\}$ are the sets from the definition of finite presentability of $G$ over $H$. The relations $s_i (a_1, \dots, a_m, X) = 1$, $i = 1, \dots, n$ are corresponded to a system of equations with coefficients from $H$:
\begin{equation*}
    \begin{cases}
        s_1 (t_1, \dots, t_m, X) = 1 \\
        \ \ \ \ \ \ \ \ \ \ \dots \\
        s_n (t_1, \dots, t_m, X) = 1 \\
    \end{cases}
\end{equation*}
which, by condition, has a solution $t_1 = a_1, \dots, t_m = a_m$. By virtue of algebraic closedness of $H$ in $G$, this system has a solution $t_1 = h_1, \dots, t_m = h_m$ in $H$. Mapping $X \sqcup \{a_1, \dots, a_m\} \to H,$ $x \in X \mapsto x$, $a_i \mapsto h_i$ extends to a surjective homomorphism $\phi: G \to H$, as defining relations of $G$ are mapped into true identities under such a mapping of generators (note that $R$ is the set of words in the alphabet $X$).

This homomorphism is the desired retraction: let $h \in H$, \ $h = v(x_1, \dots, x_r)$, $x_i \in X$. Applying to this word the homomorphism $\phi$, we get: $\phi(h) = v(\phi(x_1), \dots, \phi(x_r)) = h$. 

Algebraic closedness of a subgroup $H$ of a group $G$ follows from retractness of $H$ in $G$ for every group $G$ \cite{myasnikov}.
\qedsymbol
\end{Proof}

\begin{AppLem} \textsl{Let $C$ be a finite elementary abelian $p$-group (where $p$ is a prime number). For any $k \in \mathbb{N}$, there exists $t \geq k$ such that the direct product $P = \mbox{\Large$\times$}_{i=1}^{t} C_i$ of copies $C_i$ of $C$ contains a subgroup $R$ invariant with respect to the diagonal action on $P$ of the endomorphism algebra $\End C$ with the following properties:}
\begin{enumerate}[topsep=0pt, label=\arabic*)]
\itemsep-0.5em
\item \textsl{$R \subseteq \bigcup \ker \rho_j$, where $\rho_j : P \to C_j$, $j = 1, \dots, t$ are the natural projections,}
\item \textsl{But $R \cdot \mbox{\Large$\times$}_{j \not\in J} C_j = P$ for any subset $J \subseteq \{1, \dots, t\}$ of cardinality $|J| = k$,}
\item \textsl{Moreover, each such $J$ is contained in a set $J' \supseteq J$ such that $P = R \times (\mbox{\Large$\times$}_{j \not\in J'} C_j)$; and there exist integers $n_{ij} \in \mathbb{Z}$ such that the projection $\pi : P \to \mbox{\Large$\times$}_{j \not\in J'} C_j$ with the kernel $R$ acts as: $C_i \ni c_i \mapsto \prod_{j \not\in J'} c_j^{n_{ij}}$, where $c_j \in C_j$ is the element corresponding to $c_i$ under the isomorphism $C_i \cong C \cong C_j$.}
\end{enumerate}
\end{AppLem}

\vspace{1mm}

\noindent The following theorem provides a generalization of the result from \cite{verbal} about the center of a finite strongly verbally closed group. The proof is also analogical to the proof of that theorem, with the exception of some nuances.

\begin{FinTh}
\textsl{Let $H$ be a strongly verbally closed group. For any finite normal subgroup $T$ of $H$, for any abelian subgroup $A$ of $T$, normal in $H$, it is true that $Z(C_T (A))$ is a direct factor of $C_T (A)$, and some complement is normal in $H$. Here $C_T (X) = C(X) \cap T$.}
\end{FinTh}

\begin{Proof}
Let $H$ be such a group, and let $L = C_{T}(A)$. It suffices, for each prime $p$, to find a homomorphism $\psi_p : L \to Z(L)$ commuting with the action $H \curvearrowright L$ by conjugations (this action is well-defined as $L \unlhd H$) and injective on the $p$-component of the center $Z_p (L)$ of $L$. Then the homomorphism $ \psi : L \to Z(L) ,$ $x \mapsto \displaystyle\prod_{p} \pi_p (\psi_p (x))$, where $\pi_p : Z(L) \to Z_p (L)$ is the projection on the $p$-component, is injective on $Z(L)$, so its kernel is the desired complement $D$ (normality of $D$ in $H$ follows from the fact that $\psi$ commutes with $H \curvearrowright L$).

Suppose that there are no such homomorphisms for some prime number $p$, i.e. every homomorphism $f : L \to Z(L)$ commuting with the action $H \curvearrowright L$ is not injective on $Z_p (L)$. Then it is not injective on the maximal elementary abelian $p$-subgroup $C \leq Z_p (L)$ (it is finite as $L$ is finite). Indeed, if $x \neq 1 \in Z_p (L)$ is an element such that $f(x) = 1$, then, raising it to the appropriate power $d$, we get $f(x^d) = 1$ and $x^d \in C$, $x^d \neq 1$. 

Choose $t$ by the approximation lemma applied to $C$ (for some $k$ to be specified later) and consider the fibered product of $t$ copies of the group $H$:
\begin{equation*}
    Q = \{(h_1, \dots, h_t) \ \vert \ h_1 L = \dots = h_t L\} \leq H^t.
\end{equation*}
First of all, let us show that the subgroup $R \leq C^t \leq Q$ from the approximation lemma is normal in $Q$. Subgroup $R$ is invariant under the diagonal action of automorphisms $\Aut C \leq \End C$. It remains to show that $Q$ acts by conjugations on $P = C^t$ diagonally. It follows from the lemma:

\begin{Lem}
\textsl{Let $G$ be a group, and $N \unlhd G$. If $x C(N) = y C(N)$ for some $x, y \in G$, then $x$ and $y$ act on $N$ (by conjugations) identically.}
\end{Lem}

\begin{Proof}
From $x C(N) = y C(N)$ it follows that for some $c \in C(N)$, $x=yc$. Then for $n \in N$, we have: 
\begin{equation*}
    x^{-1} n x = c^{-1} (y^{-1} n y) c = y^{-1} n y.
\end{equation*}
The last identity is true, as (due to normality) $y^{-1} n y \in N$ and $c \in C(N)$. \qedsymbol
\end{Proof}

Let $q = (q_1, \dots, q_t) \in Q$, $p = (p_1, \dots, p_t) \in P$. As $q_1 L =$ $q_2 L =$ $\dots = q_t L$, then (according to Lemma $4$) $q^{-1} p q =$ $\Tilde{q}^{-1} p \Tilde{q}$, where $\Tilde{q} = (q_1, \dots, q_1)$. It means that the conjugation action of $Q$ on $P$ is diagonal. On the other hand, diagonal action by conjugations induces an endomorphism of $C^t$ (due to normality of $C \unlhd H$), and $R$ is invariant with respect to the diagonal action of such endomorphisms, leading to normality of $R$ in $Q$.

Put $G = Q / R$. First, let us show that $H$ embeds into $G$. The group $H$ embeds into $Q$ diagonally: $h \mapsto (h, \dots, h)$, $h \in H$. This homomorphism serves as embedding into $G$ as well, as all projections of any nontrivial diagonal element of $Q$ are nontrivial (and $R$ is contained in the union of the kernels of these projections).

Now, let us prove the verbal closedness of this diagonal subgroup (denote it as $H$ too) in $G$. Consider an equation
\begin{equation*}
    w(x_1, \dots, x_n) = (h, \dots, h)
\end{equation*}
having a solution in $G$ and let $\Tilde{x}_1, \dots, \Tilde{x}_n$ be a preimage (in $Q$) of a solution $x_1, \dots, x_n$. Then (in $Q$):
\begin{equation*}
    w(\Tilde{x}_1, \dots, \Tilde{x}_n) = (hc_1, \dots, hc_t),
\end{equation*}
where $(c_1, \dots, c_t) \in R$. By the property $1)$ of the approximation lemma, $c_i = 1$ for some $i$. It means that in $H$ (the group itself) $w(\Tilde{x}_1^i, \dots, \Tilde{x}_n^i) = h$, where $\Tilde{x}_j^i$ is the $i$th coordinate of the vector $\Tilde{x}_j$, $j = 1, \dots, n$.

Let us take $y_j = (\Tilde{x}_j^i, \dots, \Tilde{x}_j^i)$, $j = 1, \dots, n$. Then in $H \leq G$ the following is true:
\begin{equation*}
    w(y_1, \dots, y_n) = (h, \dots, h),
\end{equation*}
which proves the verbal closedness of $H$ in $G$.

Let $U \leq L$. We use the following denotion:
\begin{equation*}
    U_i : = \{(1, \dots, 1, u, 1, \dots, 1) \ \vert \ u \in U \} \leq Q, \ i = 1, \dots, t \text{ (coordinate } u \text{ stands on the } i \text{th place)}.
\end{equation*}
It remains to prove that $H$ is not algebraically closed in $G$. 

\begin{Lem}
\textsl{The group $Q$ is finitely presented over its subgroup $H$.}
\end{Lem}

\begin{Proof} According to Lemma $2$, it is sufficient to show that $Q / (L_1 \times \dots \times L_t)$ is finitely presented over $H / \Tilde{L}$, where $\Tilde{L} = \{(l, \dots, l) \ \vert \ l \in L\}$. However, $Q = H \cdot (L_1 \times \dots \times L_t)$, so the statement we prove follows from this fact (see \cite{kargapolov}, theorem $4.2.4$): 
\begin{quote}
    Suppose that $G$ is a group, $F$ is its subgroup, and $K$ is its normal subgroup. Then $(K \cdot F) / K \cong F / (F \cap K)$.
\end{quote}
Thus, the group $Q / (L_1 \times \dots \times L_t)$ is not just finitely presented over $H / \Tilde{L}$ but is isomorphic to it. \qedsymbol
\end{Proof}

From Lemma $3$ and Lemma $5$, it follows that it suffices to show that $H$ is not a retract of $G$. Let $\rho : G \to H$ be a hypothetical retraction, and let $\hat{\rho} : Q \to H$ be its composition with the natural epimorphism $Q \to Q / R = G$. Henceforth, all subgroups and centralizers we refer to relate to $Q$. 

Let us verify that $\hat{\rho} (L_i) \leq C_{T} (C_{T} (L)) \leq L$ for every $i$. First, prove the left inclusion. Let $h \in C_{T} (L)$. Then, $h$ commutes with every element from $L$; consequently, $h$, as an element of $Q$, commutes with $L_i$. Applying the retraction $\hat{\rho}$ to this identity, we get that $\hat{\rho} (h)$ $(= h)$ commutes with the subgroup $\hat{\rho} (L_i)$, which (by definition of the centralizer) proves the inclusion. The second inclusion follows from the fact that $L = C_{T}(A) =$ $C(A) \cap T$, which means that
\begin{equation*}
    C_{T} (C_{T} (L)) \leq C_{T} (A \cap T) = C_{T} (A) = L.
\end{equation*}
The first inclusion here is true as $C(L) \geq A$. The following equality is true as $A \leq T$.

On the other hand, for $i \neq j$, the mutual commutator subgroup $[L_i, L_j]$ is trivial (as in case $i$ and $j$ are different, $L_i$ and $L_j$ are contained in different components of the fibered product). It means that the image of this mutual commutator subgroup is trivial too: $[\hat{\rho} (L_i), \hat{\rho} (L_j)] = \{1\}$. Consequently, $[L_i, \prod_{j \neq i} L_j] = \{1\}$ and $[\hat{\rho}(L_i), \prod_{j \neq i} \hat{\rho}(L_j)] = \{1\}$. If $\hat{\rho}(L_i) = \hat{\rho}(L_l)$ for some $i \neq l$, then (by the virtue of well-known commutator identities) $[\hat{\rho}(L_i), \prod_{j} \hat{\rho}(L_j)] = \{1\}$, which means that $\hat{\rho} (L_i) \leq C_{T} (L)$ (as $L = \hat{\rho} (L) \leq \prod_{j} \hat{\rho} (L_j)$).

Thereby, if for some different $i$ and $j$, $\hat{\rho}(L_i) = \hat{\rho}(L_j)$, then $\hat{\rho} (L_i) \leq C_{T} (L)$. From here and from the inclusion we proved earlier, we get $\hat{\rho} (L_i) \leq L \cap C_{T}(L) =$ $Z (L)$.

Let us take $k$ in the approximation lemma to be the number of all subgroups of $T$, and let $J$ be the set of all \textit{exclusive} numbers $i$, namely such that for any $l \neq i$, $\hat{\rho} (L_i) \neq \hat{\rho} (L_l)$. Since among $\hat{\rho} (L_i) \leq T$ there are no more than $k$ different subgroups, $|J| \leq k$. Thus, from the property $3)$ of the approximation lemma, we have a decomposition:
\begin{equation*}
    \mbox{\Large$\times$}_{i = 1}^t C_i = R \times (\mbox{\Large$\times$}_{i \in I} C_i),
\end{equation*}
where $I \subseteq \{1, \dots, t\} \setminus J$ is some set of non-exclusive elements. Again, according to the property $3)$ of the approximation lemma, the projection $\pi : \mbox{\Large$\times$}_{i = 1}^t C_i \to \mbox{\Large$\times$}_{i \in I} C_i$ onto the second factor of this decomposition is defined by an integer matrix $(n_{ij})$, namely, for $c_i \in C_i$, $\pi : c_i \mapsto \prod_{j \in I} c_j^{n_{ij}}$, where $c_j$ are elements corresponding to $c_i$ under the isomorphism $C_i \cong C \cong C_j$.

This means that the restriction of $\pi$ to $C = \{(c, \dots, c) \ \vert \ c \in C \leq H\}$ is defined by formula:
\begin{equation*}
    \hat{\pi} : (c, \dots, c) \mapsto \prod_{j \in I} c_j^{m_j}, \ m_j = \sum_{i} n_{ij}.
\end{equation*}
Here $c_j$ are elements corresponding to $c$ under the isomorphism $C \cong C_j$.

Then (as $i \in I$ are non-exclusive, we have $\hat{\rho} (L_i) \leq Z(L)$), consider the composition:
\begin{equation*}
    \Psi : C \leq Q \to Z(L), \ c \overset{\pi}{\mapsto} \prod_{j \in I} c_j^{m_j} \overset{\hat{\rho}}{\mapsto} \prod_{j \in I} \hat{\rho}(c_j^{m_j}).
\end{equation*}
It extends to a homomorphism $\Phi : L \to Z(L)$ defined by the similar formula:
\begin{equation*}
    \Phi : g \mapsto \prod_{j \in I} \hat{\rho}(g_j^{m_j}),
\end{equation*}
where $g \in L$ and $g_j \in L_j$ are elements corresponding to $g \in L$. Obviously, it is an extension of $\Psi$ and a homomorphism, as for $j \in I$, $\hat{\rho} (L_j) \leq Z(L)$ and the group $Z(L)$ is abelian. This homomorphism commutes with the conjugation action of $H$ on $L$. Indeed, let $g \in H$ and let $\mathfrak{g}$ be the action of $g$ on $L$ by conjugation, namely, for $x \in L$, $\mathfrak{g} (x) = g^{-1} x g$. Let us show that $\Phi \circ \mathfrak{g} = \mathfrak{g} \circ \Phi$. Let $h \in L$. Then $\Phi (\mathfrak{g}(h)) =$ $\prod_{j \in I} \hat{\rho}(g^{-1} h_j^{m_j} g) =$ $\prod_{j \in I} g^{-1} \hat{\rho}(h_j^{m_j}) g =$ $\mathfrak{g} (\Phi (h))$. Penult identity is true, as $\hat{\rho}$ is a retraction on $H$, so it acts identically on $H$ itself. By assumption we made in the beginning, the kernel of this homomorphism has nontrivial intersection with $C$: $\ker \Phi \cap C \neq \{1\}$, so the restriction $\Psi = \Phi \vert_C$ has a nontrivial kernel too.

On the other hand, $\Psi$ is the identical mapping, since $\Psi = \hat{\rho} \vert_C \circ \pi \vert_C = $ $\hat{\rho} \vert_C \circ \hat{\pi} = $ $\hat{\rho} \vert_C$ (the last identity is true as $\hat{\pi}$ is a projection <<forgetting>> the $R$ coordinate, and $\hat{\rho} (R) = \{1\}$ is a composition of the natural homomorphism to the quotient group and of the retraction to $H$) and $\hat{\rho} \vert_C = id$, as $\hat{\rho}$ is the retraction from $Q$ to $H$, so it acts trivially on $C$. The obtained contradiction completes the proof. \qedsymbol

\end{Proof}

\noindent Let us provide some corollaries of this theorem:

\begin{Cons}
\textsl{Finitely generated nilpotent groups with nonabelian torsion subgroups are not strongly verbally closed.}
\end{Cons}

\begin{Proof}
Let us take the torsion subgroup of such group as $T$ from the theorem, and the center of this torsion subgroup as $A \neq \{1\}$. Since $T$ is nilpotent and nonabelian, every nontrivial normal subgroup of $T$ has a nontrivial intersection with $A$ \cite{kargapolov}, so $A$ is not a direct factor of $T$. \qedsymbol
\end{Proof}

\begin{Cons2}
\textsl{A finite group, whose center is not a direct factor is not strongly verbally closed.}
\end{Cons2}

\vspace{1mm}

This theorem does not cover the case of finitely generated nilpotent nonabelian groups with abelian torsion subgroups, and it is still unknown whether there are strongly verbally closed groups among such groups. So far, we can provide only a partial answer to this question (see the first proposition of the following paragraph).

\section*{\normalsize 4. Nilpotent non-strongly-verbally-closed groups}

\noindent Let us remind that \textit{the discrete Heisenberg group} is the free nilpotent group of the nilpotency class two with two free generators. It can be easily verified that this group admits a faithful representation in the group of upper triangular matrices of size $3$ by $3$.
\begin{Prop}
\textsl{Let $H$ be the discrete Heisenberg group with $a$ and $b$ being its free generators and $N$ being its subgroup:}
\begin{equation*}
    N = \langle \langle a^{\alpha}, [a,b]^n \rangle \rangle, \ \alpha, n \geq 0.
\end{equation*}
\noindent \textsl{The group $G = H / N$ is strongly verbally closed if and only if} $\text{gcd}(\alpha, n) = 1$.
\end{Prop}

\begin{Proof}
Let $T(G)$ be the torsion subgroup of $G$. The center of the group $H$ is equal to its commutator subgroup, and it is isomorphic to the infinite cyclic group. As it was said earlier (refer to the proof of the Proposition $4$), any nontrivial normal subgroup of a nilpotent group intersects its center nontrivially. Thus, if $T(G) = \{1\}$, then either $G = H$ is the discrete Heisenberg group, or $G$ is abelian. Non-strong-verbal-closedness of $H$ was proved in \cite{verbal}, and the strong verbal closedness of abelian groups  was proved in \cite{mazhuga}. The case of $T(G) = \{1\}$ corresponds to the case of $\alpha = 0$, $n = 0$ and $\alpha = 0$, $n = 1$. 

If $\text{gcd}(\alpha, n) = 1$, then, once again, $G$ is abelian, since $[a,b]^{\alpha} = [a^{\alpha}, b] \in N$; consequently, it is strongly verbally closed. 

Consider the case, when $\text{gcd}(\alpha, n) = d \neq 1$. Without loss of generality, we may assume that $\alpha$ and $n$ are the least of non-negative numbers such that $a^{\alpha} \in N$, $[a,b]^{n} \in N$. Consider the central product of $G$ with its copy $\Tilde{G}$ with joined commutator subgroup:
\begin{equation*}
    K = G \underset{G' = \Tilde{G}'}{\times} \Tilde{G}= (G \times \Tilde{G}) / \{(c, c^{-1}) \vert c \in G'\}.
\end{equation*}
The group $G$ is not algebraically closed in $K$, since $G$ is not a retract of $K$. Indeed, let $\rho$ be a hypothetical retraction. The group $G$ commutes with $\Tilde{G}$ in $K$, so $\rho(\Tilde{G}) \leq Z(G)$ and $\rho(\Tilde{G}') = \{1\}$, which leads to a contradiction with the definition of retraction. However, $G$ is verbally closed in $K$. Let $w \in F(t_1, \dots, t_s)$ be some word and
\begin{equation*}
    w((h_1 N, h_1' N), \dots, (h_s N, h_s' N)) = (h N,N)
\end{equation*}
for some $h N, h_i N \in G$, $h_i' N \in \Tilde{G}$. Then, for some $cN \in G'$, the following holds:
\begin{equation*}
    \begin{cases}
    w(h_1', \dots, h_s' ) N = cN \\
    w(h_1, \dots, h_s) N = hc^{-1} N \\
    \end{cases}
\end{equation*}
By an automorphism of the free group, the word $w$ can be reduced to a \textit{normal form} \cite{verbal}: $w(t_1, \dots, t_s) = t_1^m w'(t_1, \dots, t_s)$, where $m \in \mathbb{N}$, $w' \in F_s'$. From the first equation, we get $c N \in G' \ \cap \ \phi(G^s)$, where $\phi : G^s \to G,$ $(g_1, \dots, g_s)$ $\mapsto w(g_1, \dots, g_s)$ is a verbal mapping. This means that for some $w_1, w_2 \in N$, in $H$ it is true that:
\begin{equation*}
    \begin{cases}
    w(h_1', \dots, h_s') = c w_1 \\
    w(h_1, \dots, h_s) = hc^{-1} w_2 \\
    \end{cases}
\end{equation*}
Let us show that in $G$ the identity $(aN)^x = [aN,bN]^{z}$ doesn't hold for $x \not\in \alpha \mathbb{Z}$. Converse would mean that in the discrete Heisenberg group the following holds:
\begin{equation*}
    a^x [a,b]^{-z} = b^{-k} a^{-l} a^{\alpha t} [a,b]^{n s} a^{l} b^{k}
\end{equation*}
for some $k, l, t, s \in \mathbb{Z}$. After some reductions, we get: $a^{x - \alpha t} = [a,b]^{ns + z + \alpha t k}$. In $H$ it is possible only if $x = \alpha t \in \alpha \mathbb{Z}$. We obtained a contradiction. 
Thus, $h_1'= [a,b]^{\gamma}$ for some $\gamma \in \mathbb{Z}$, and, consequently, $c w_1 \in H'$. Since for any verbal mapping $\phi$ in the discrete Heisenberg group (see \cite{verbal}) 
\begin{equation*}
    \text{for any } g \in \phi (H^s),\text{ it is true that } g (\phi(H^s) \cap H') \subseteq \phi(H^s),
\end{equation*}
for some $g_1, \dots, g_s \in H$, $w(g_1, \dots, g_s) =$ $w(h_1, \dots, h_s) c w_1$. It means that:
\begin{equation*}
    w(g_1, \dots, g_s) = h w_3
\end{equation*}
for some $w_3 \in W$, and in $G$:
\begin{equation*}
    w(g_1, \dots, g_s) N = h N,
\end{equation*}
which proves verbal closedness of $G$ in $K$. \qedsymbol
\end{Proof}

\vspace{1mm}

At last, let us prove that the higher dimensional Heisenberg groups over any field are not strongly verbally closed:

\textit{The Heisenberg group of dimension $2n+1$ over a field $K$}, where $n \in \mathbb{N}$ is the group of upper triangular matrices of the kind
\begin{equation*}
H_n (K) =
\Bigg\{ T(\Bar{a}, \Bar{b}, c) =
    \begin{pmatrix}
    1 & \Bar{a} & c \\
    0 & I_n & \Bar{b} \\
    0 & 0 & 1 \\
    \end{pmatrix} \Bigg| \Bar{a}, (\Bar{b})^{\intercal} \in K^n, c \in K
\Bigg\},
\end{equation*}
where $I_n$ is the identity matrix of size $n$.

\vspace{1mm}

\begin{Prop}
\textsl{The group $H_n (K)$ is not strongly verbally closed.}
\end{Prop}

\begin{Proof}
Consider the central product of $H_n (K)$ with its copy $\Tilde{H}_n (K)$ with joined commutator subgroup:
\begin{equation*}
    G = H_n (K) \underset{H_n (K)' = \Tilde{H}_n (K)'}{\times} \Tilde{H}_n(K).
\end{equation*}

\noindent Denote with the symbol $H$ the first factor of this central product. Let us show that $H$ is not algebraically closed in $G$. The group $H$ is linear, and, consequently, it is \textit{equationally noetherian} \cite{baumslag}, so it is algebraically closed in $G$ if and only if it is a retract of every finitely generated over $H$ subgroup of $G$ \cite{virtually}. In particular, of such a subgroup of $G$:
\begin{equation*}
\Bar{H} = \langle H, (1, h_1), \dots, (1, h_n), (1, g_1), \dots, (1, g_n) \rangle \text{, where }
    h_i = \begin{pmatrix}
    1 & \Bar{a}_i & 0 \\
    0 & I_n & 0 \\
    0 & 0 & 1 \\
    \end{pmatrix}, g_i =
    \begin{pmatrix}
    1 & 0 & 0 \\
    0 & I_n & \Bar{b}_i \\
    0 & 0 & 1 \\
    \end{pmatrix},
\end{equation*}
where $\Bar{a}_i = (0, \dots, 1, \dots, 0) =(\Bar{b}_i)^{\intercal}$, (unit is on the $i$th place). Thus, $N = \langle h_1, \dots, h_n, g_1, \dots, g_n \rangle$ is a subgroup of $\tilde{H}_n (K)$, isomorphic to the discrete Heisenberg group of the dimension $(2n + 1)$. Let $\rho$ be a hypothetical retraction. Since in $G$ the group $H$ commutes with $N$, we get that $\rho(N') = \{1\}$, which leads to a contradiction with the definition of retraction.

Nevertheless, subgroup $H$ is verbally closed in $G$: let $w \in F_s$ be some word (without loss of generality, this word is in the normal form we established earlier), and let $\phi : H^s \to H$ be the verbal mapping associated with this word. Suppose that for some $h_i, h \in H$, $h_i' \in \Tilde{H}$, $c \in H'$:
\begin{equation*}
    \begin{cases}
    w(h_1', \dots, h_s') = c \\
    w(h_1, \dots, h_s) = hc^{-1} \\
    \end{cases}
\end{equation*}
In general, on matrices $g_i = T(\bar{a}_i, \bar{b}_i, c_i)$, $i=1, \dots, s$, the mapping $\phi$ acts like that:

\begin{equation*}
    \phi(g_1, \dots, g_s) = \begin{pmatrix}
        1 & m \bar{a}_1 & m c_1 + f(\bar{a}_1, \dots, \bar{a}_s; \bar{b}_1, \dots, \bar{b}_s) \\
        0 & I_n & m \bar{b}_1\\
        0 & 0 & 1 \\
    \end{pmatrix},
\end{equation*}
where $f: (K^n)^s \times (K^n)^s \to K$ is some function linear in every argument. The image of $f$ is either trivial or is equal to $K$, which leads to:

\begin{equation*}
    \phi(H_n (K)^s) = \begin{cases}
    \{1\}, & \text{ if } m = 0 \text{ and the image of } f \text{ is trivial} \\
    (H_n (K))',  & \text{ if } m = 0 \text{ and the image of } f \text{ equals } K \\
    H_n (K),  & \text{ if } m \neq 0
    \end{cases}
\end{equation*}
Then $\phi(H_n (K)^s) \cap (H_n (K))' \leq H_n (K)$ and for every element $h \in \phi(H_n (K)^s)$ it is true that
\begin{equation*}
    h(\phi(H_n (K)^s) \cap (H_n (K))') \subseteq \phi(H_n (K)^s),
\end{equation*}
whence verbal closedness follows.
\qedsymbol

\end{Proof}


\begingroup
\setlength\bibitemsep{0pt}
\printbibliography[title={\large{\textmd{REFERENCES}}}]
\endgroup

\end{document}